# An Introduction To The Theory Of Newton Polygons For L-Functions Of Exponential Sums

Daqing Wan

Department of Mathematics
University of California
Irvine, CA 92697-3875
dwan@math.uci.edu

**Table of Contents**

This paper is based on the author's series of lectures delivered at the January 1999 Mini-course in Number Theory, held at Sogang University (Seoul). The aim was to give an elementary and self-contained introduction to the theory of Newton polygons (namely, the $p$-adic Riemann hypothesis) for L-functions of exponential sums over a finite field. In addition to giving a thorough treatment of the basic elementary local theory, we also describe the deeper global theory and include some explicit examples to illustrate how to use the main theorems. In this way, it is hoped that these notes would at least give a little feeling and some indications about the subject of Newton polygons for L-functions. This arithmetic subject can be traced back to Stickelberger's classical theorem [St] on the prime factorization of a Gauss sum. Since then, tremendous developments have been made in many directions by many people. Despite these progresses, there are still many more fundamental open questions than theorems. The situation is unlikely to be drastically improved anytime soon due to the diversity, richness and depth of both the questions and the methods involved in this subject. To some extent, it is also due to our general lack of sufficient experiences with these subtle arithmetic problems.

The main idea of our approach is to establish a suitable local to global principal which would allow us to reduce the harder "global" case to various easier well understood "local" cases for which Stickelberger's theorem applies. For this purpose, a systematic decomposition method is introduced in [W2-3] and several decomposition theorems are proved, including the facial decomposition theorem, the star decomposition theorem and the hyperplane decomposition theorem. These theorems easily recover previously known results on Newton polygons of exponential sums. They are also enough for a number of further applications such as Mazur's conjecture [Ma1] for a generic hypersurface and a weaker form of the more general Adolphson-Sperber conjecture [AS3] for a generic exponential sum. Using the method in [W3], we have recently found (in fact, during the initial preparation of these lecture notes) a more flexible collapsing decomposition theorem. This theorem gives a fairly satisfactory answer to the full form of the Adolphson-Sperber conjecture. In particular, as an application, we now know that the full form of the Adolphson-Sperber conjecture is true in every low dimension $n \leq 3$ but false in every high dimension $n \geq 4$. These results are explained in these notes.

## Chapter 1. Introduction

In this introductory chapter, we first review the basic theory of L-functions of exponential sums over finite fields and our main question on the Newton polygon problem. Then, we turn to a general intuitive discussion about the available results and some open questions. This part can be viewed as a quick exposition of this subject without a detailed description.

### 1.1. L-functions of exponential sums over finite fields.

Let $\mathbf{F}_q$ be the finite field of $q$ elements with characteristic $p$. For each positive integer $k$, let $\mathbf{F}_{q^k}$ be the finite extension of $\mathbf{F}_q$ of degree $k$. Let $\zeta_p$ be a fixed primitive $p$-th root of unity in the complex numbers. For any Laurent polynomial $f(x_1, \cdots, x_n) \in \mathbf{F}_q[x_1, x_1^{-1}, \cdots, x_n, x_n^{-1}]$, we form the exponential sum

$$S_k^*(f) = \sum_{x_i \in \mathbf{F}_{q^k}^*} \zeta_p^{\mathrm{Tr}_k f(x_1, \cdots, x_n)}, \tag{1.1}$$

where $\mathbf{F}_{q^k}^*$ denotes the set of non-zero elements in $\mathbf{F}_{q^k}$ and $\mathrm{Tr}_k$ denotes the trace map from $\mathbf{F}_{q^k}$ to the prime field $\mathbf{F}_p$. This is an exponential sum over the $n$-torus $\mathbf{G}_m^n$ over $\mathbf{F}_{q^k}$. A fundamental question in number theory is to understand the sequence $S_k^*(f)$ ($1 \leq k < \infty$) of algebraic integers, each of them lying in the $p$-th cyclotomic field $\mathbf{Q}(\zeta_p)$.

By well known theorems of Dwork-Bombieri-Grothendieck, the following generating L-function is a rational function:

$$L^*(f,t) = \exp\Big(\sum_{k=1}^{\infty} S_k^*(f)\frac{t^k}{k}\Big) = \frac{\prod_{i=1}^{d_1}(1-\alpha_i t)}{\prod_{j=1}^{d_2}(1-\beta_j t)}, \tag{1.2}$$

where the finitely many numbers $\alpha_i$ $(1 \leq i \leq d_1)$ and $\beta_j$ $(1 \leq j \leq d_2)$ are non-zero algebraic integers. Equivalently, for each positive integer $k$, we have the formula

$$S_k^*(f) = \beta_1^k + \beta_2^k + \cdots + \beta_{d_2}^k - \alpha_1^k - \alpha_2^k - \cdots - \alpha_{d_1}^k.$$

Thus, our fundamental question is reduced to understand the reciprocal zeros $\alpha_i$ $(1 \leq i \leq d_1)$ and the reciprocal poles $\beta_j$ $(1 \leq j \leq d_2)$. When we need to indicate the dependence of the L-function on the ground field $\mathbf{F}_q$, we will write $L^*(f/\mathbf{F}_q, t)$.

Without any condition on $f$, one does not even know exactly the number $d_1$ of zeros and the number $d_2$ of poles, although good upper bounds are available, see [Bom] and [AS2]. On the other hand, Deligne's well known theorem [De2] gives the following general information about the nature of the zeros and poles. For the complex absolute value $|\ |$, this says

$$|\alpha_i| = q^{u_i/2}, \ |\beta_j| = q^{v_j/2}, \ u_i \in \mathbf{Z} \cap [0, 2n], \ v_j \in \mathbf{Z} \cap [0, 2n], \tag{1.3}$$

where $\mathbf{Z} \cap [0, 2n]$ denotes the set of integers in the interval $[0, 2n]$. For each $\ell$-adic absolute value $|\ |_\ell$ with prime $\ell \neq p$, the $\alpha_i$ and the $\beta_j$ are $\ell$-adic units:

$$|\alpha_i|_\ell = |\beta_j|_\ell = 1. \tag{1.4}$$

For the remaining prime $p$, (1.3) and (1.4) show

$$|\alpha_i|_p = q^{-r_i}, \ |\beta_j|_p = q^{-s_j}, \ r_i \in \mathbf{Q} \cap [0, n], \ s_j \in \mathbf{Q} \cap [0, n], \tag{1.5}$$

where we have normalized the $p$-adic absolute value by $|q|_p = q^{-1}$. Strictly speaking, in defining the $p$-adic absolute value, we have tacitly chosen an embedding of the field $\bar{\mathbf{Q}}$ of algebraic numbers into an algebraic closure of the $p$-adic number field $\mathbf{Q}_p$. The slightly weaker version of (1.5) with $[0, n]$ replaced by $[0, 2n]$ also follows easily from a trivial archimedian estimate. The precise version of various types of Riemann hypothesis for the L-function in (1.2) is then to determine the important arithmetic invariants $\{u_i, v_j, r_i, s_j\}$. The integer $u_i$ (resp. $v_j$) is called the weight of the algebraic integer $\alpha_i$ (resp. $\beta_j$). The rational number $r_i$ (resp. $s_j$) is called the slope of the algebraic integer $\alpha_i$ (resp. $\beta_j$) defined with respect to $q$. Without any condition on $f$, not much more is known about these weights and the slopes, since one does not even know exactly the number $d_1$ of zeros and the number $d_2$ of poles.

Assume now that $f$ is non-degenerate (a smooth condition on the "leading form" of $f$) with respect to $\Delta(f)$, where $\Delta(f)$ is the $n$-dimensional convex polyhedron in $\mathbf{R}^n$ generated by the origin and the exponents of the non-zero monomials which occur in $f$. Then, the L-function $L^*(f,t)^{(-1)^{n-1}}$ is known to be a polynomial of degree $n!\mathbf{V}(f)$ by a theorem of Adolphson-Sperber [AS3] proved using $p$-adic methods, where $\mathbf{V}(f)$ denotes the volume of $\Delta(f)$. The complex absolute values (or the weights) of the $n!\mathbf{V}(f)$ zeros can be determined explicitly by a theorem of Denef-Loeser [DL] proved using $\ell$-adic methods. As indicated above, the $\ell$-adic absolute values of the zeros are always 1 for each prime $\ell \neq p$. Thus, there remains the intriguing question of determining the $p$-adic absolute values (or the slopes) of the zeros. This is the $p$-adic Riemann hypothesis for the L-function $\mathrm{L}^*(f,t)^{(-1)^{n-1}}$. Equivalently, the question is to determine the Newton polygon of the polynomial $\mathrm{L}^*(f,t)^{(-1)^{n-1}}$. This $p$-adic question turns out to be extremely complicated in general, although a standard lower bound of the Newton polygon is known [AS3], in the direction of Katz's well known conjecture [Ka1] that the Newton polygon lies above the Hodge polygon. The necessary complication of the Newton polygon problem is partly due to the fact that there does not exist a clean general answer or equivalently the answer varies too much as $f$ and the prime $p$ vary, although Grothendieck's specialization theorem does provide a general structure theorem about the variation of the Newton polygon when $f$ varies in an algebraic family. This suggests that it may be necessary to restrict to certain classes of polynomials $f$ (hopefully, of reasonable generality) and certain classes of $p$ in order to obtain a clean and precise result about the Newton polygon.

### 1.2. A brief history on the Newton polygon problem.

In the simplest case when $f$ is the one variable monomial $x^d$, the exponential sum in (1.1), sometimes also called the Gauss sum in the literature, can be easily expressed in terms of the standard Gauss sums involving both additive and multiplicative characters of $\mathbf{F}_q$, see Berndt and Evans [BE] for a comprehensive survey and [BEW] for a complete updated account on this classical topic. The rationality of the L-function in this case follows from the Hasse-Davenport relation [HD] on Gauss sums. It shows that each reciprocal zero of the L-function in this case is a radical of a Gauss sum. The $p$-adic absolute value of a Gauss sum (and hence its radical) is determined by Stickelberger's theorem [St]. The final result in this direction is the precise $p$-adic formula of Gross-Koblitz [GK] for the Gauss sum in terms of the $p$-adic $\Gamma$-function.

More generally, if $f$ is a higher dimensional diagonal Laurent polynomial, the L-function can also be computed explicitly using Gauss sums and the Hasse-Davenport relation. This generalizes the classical case of the zeta function of a Fermat hypersurface as calculated by Weil [We] to support his celebrated conjectures. Each reciprocal zero of the L-function in this case is a radical of a product of Gauss sums. Combining this with the Stickelberger theorem, one obtains, at least theoretically, a complete understanding of the Newton polygon problem in the elementary diagonal case. A thorough treatment of the diagonal case is given in Chapter 3, which will play the role of building blocks (local cases) for later more general global theory. Despite the elementary nature of the diagonal case, there may still be something to be learned from it. In particular, non-trivial arithmetic and combinatorial problems often arise in the actual explicit calculation of the Newton polygon even in the diagonal case due to the diversity of $\Delta(f)$ and the nature of the prime $p$. In addition, the diagonal case provides an excellent (checkable) testing ground and a rich source of explicit examples for conjectures on L-functions.

For non-diagonal $f$, the Newton polygon problem is significantly deeper because there are no known explicit formulas for the reciprocal zeros. The geometric elliptic curve case was handled by Deuring [D]. For exponential sums, the first non-diagonal example seems to be the $n$-dimensional Kloosterman polynomial

$$f = a_1 x_1 + \cdots + a_n x_n + \frac{b}{x_1 \cdots x_n}$$

as completed by Dwork [Dw3] for $n = 1$ and by Sperber [S1] for all $n$. Further extensive and in depth investigations of the Newton polygon problem for various generalized Kloosterman polynomials were carried out by Sperber [S2] and Carpentier [Ca]. These special but precise results led eventually to Adolphson-Sperber's fairly general conjecture [AS3] about the Newton polygon of a generic exponential sum. This conjecture suggests a deeper relationship between the arithmetic Newton polygon and the geometric (combinatorial) Hodge polygon. It includes Mazur's earlier conjecture [Ma1] on the zeta function of a generic hypersurface as a special case.

### 1.3. Several decomposition theorems.

In [W2], a facial decomposition theorem for the Newton polygon is proved. This theorem is in some sense a local to global principal. It allows us to cut $f$ into simpler pieces if $\Delta(f)$ has more than one co-dimension 1 faces not containing the origin. If each of the simpler pieces is diagonal, then we get non-trivial information about the Newton polygon for the original global $f$. Applying this decomposition to Kloosterman type polynomials, one immediately recovers the results of Dwork-Sperber-Carpentier. It also shows that one can sometimes deform $f$ a little bit without changing the Newton polygon. As an application, one can construct many new classes of non-diagonal Laurent polynomials for which the Newton polygon of the L-function can be determined explicitly, including the new bi-Kloosterman polynomial arising from the recent calculation by Kim [Ki] of generalized Gauss sums for unitary groups. These results are discussed in Section 4.1.

In addition to treating a single Laurent polynomial $f$ with a given $\Delta(f)$, we can also consider the universal family of non-degenerate Laurent polynomials $f$ with a fixed $n$-dimensional polyhedron $\Delta(f) = \Delta$. In this family case, one can show that the L-function $L^*(f,t)^{(-1)^{n-1}}$ has a generic Newton polygon by Adolphson-Sperber's result and Grothendieck's specialization theorem for F-crystals [Ka2]. The conjecture of Adolphson-Sperber [AS3] says that if $p$ is a prime number in the residue class $p \equiv 1 \pmod{D}$ for a certain explicit $D$ depending in a simple way on $\Delta$, then the generic Newton polygon coincides with its lower bound (the so-called ordinary case). The denominator $D$ is in general optimal and cannot be improved. For instance, in all geometric situations (the zeta function case), the denominator $D$ is 1. Thus, the Adolphson-Sperber conjecture includes Mazur's earlier conjecture [Ma1] for the zeta function of a generic hypersurface as a special case.

In order to show that the generic Newton polygon coincides with its lower bound, we need to show the non-vanishing of the various Hasse polynomials which give the relevant $p$-adic information about the L-function. These Hasse polynomials in many variables are quite complicated in general. A natural idea is to show that certain simpler leading form of such a Hasse polynomial is non-zero. For this purpose, we introduced the finer notion of degree polygon in [W3] constructed using the degree of the leading form defined in terms of certain parameters, where the degree may be $-\infty$ if the leading form vanishes. An upper bound is given for the degree polygon. Although the Newton polygon goes up under specialization, the degree polygon goes down under specialization. The generic Newton polygon will coincide with its lower bound whenever the degree polygon coincides with its upper bound.

We are thus led to study when the degree polygon coincides with its upper bound. Again, the basic idea is to establish a suitable local to global principal for the degree polygon. In [W3], two decomposition theorems (the star decomposition and the hyperplane decomposition) for the degree polygon are proved. Repeatly applying the star decomposition, one can then cut a generic $f$ into smaller and smaller pieces, eventually into diagonal cases. In this way, we obtain a systematic method to determine the degree polygon and hence the generic Newton polygon. In particular, it implies that a weaker version of the Adolphson-Sperber conjecture is true. Using the method in [W3], we have recently proved a more flexible collapsing decomposition theorem for the degree polygon. This new theorem gives significantly better results in studying the full form of the Adolphson-Sperber conjecture. The detail proof will be included elsewhere. In Section 4.2, we explain the new collapsing decomposition and some of its applications. In Section 4.3, we review the old hyperplane decomposition [W3] and its application to the case of Mazur's conjecture.

The new collapsing decomposition theorem also implies the weaker version of the Adolphson-Sperber conjecture which says that the conjecture is true for all $p$ in the residue class $p \equiv 1 \pmod{D^*}$ for some finite integer $D^* \geq D$. The smallest possible $D^*$ is much more subtle than the explicit $D$. It depends a lot on how the prime $p$ relates to the arithmetic properties of the rational points in $\Delta$. Our method only produces a good effective upper bound for the optimal $D^*$. Instead of considering the largest universal family, we could consider the sub-family of the non-degenerate $f$ with a given $\Delta$ but some designated terms of $f$ being omitted. Our decomposition theorems are proved for such a more general sub-family. The integer $D$ will be the same but the desired quantity $D^*$ may increase because the indecomposable local pieces may have changed.

The full form of the Adolphson-Sperber conjecture is in general false for each $n \geq 4$, already in certain indecomposable diagonal cases. Counter-examples are given in Section 3.4. The full form is nevertheless true in many important cases such as in every low dimension $n \leq 3$ and in the geometric case of Mazur's conjecture. Usually, our decomposition method proves the full form of the Adolphson-Sperber conjecture whenever it is true for a given $\Delta$. For an indecomposable simplex, the lattice points in the fundamental domain forms a finite abelian group. The full conjecture is true for the indecomposable simplex if the largest invariant factor of the associated finite abelian group is a factor of $D$. For an arbitrary $\Delta$, the full conjecture will be true if, in our complete collapsing decomposition of $\Delta$ into indecomposable simplicial ones, the largest invariant factor associated to each simplex is a factor of $D$. This last property also occurs in the resolution of singularities of a toric variety [KK] in which case $D = 1$. The precise general relation in this connection does not seem to have been exploited.

### 1.4. Further remarks on open questions.

Our basic setting in this paper is the L-function of the exponential sum of a Laurent polynomial over the $n$-torus. We choose to work over the $n$-torus here because the situation is simplest to describe and because many other cases can often be reduced to the $n$-torus case, see [W3]. In addition to the exponential sum case discussed here, we note that there are many other results and methods available (particularly for the more elementary slope zero part of the Newton polygon) in various geometric situations such as curves, $K$-3 surfaces, abelian varieties, complete intersections, see Artin [Ar], Illusie [Il], Koblitz [Ko], Manin [M], Miller [Mi], Mumford [Mu], Norman-Oort [NO], Ogus [Og1], Oort [Or]. In particular, the first proof of Mazur's conjecture for a generic projective complete intersection was given by Illusie using some of Deligne's ideas [De2] and the DeRham-Witt complex. Our method, based on Dwork's trace formula and Adolphson-Sperber's theorem, is completely different. It works for open varieties and exponential sums as well.

It should be noted that good lower bounds for the first non-trivial slope of the Newton polygon, such as the Ax-Katz type theorem [Ax][Ka1] and its generalizations by Adolphson-Sperber [AS1] and Moreno-Moreno [MM], can all be treated in elementary way [W1] by Gauss sums and the Stickelberger theorem, or even simpler, by a version of the classical Chevalley-Warning type argument, see [W4] and the references listed there. The author believes that there is still potential to obtain interesting further results of reasonable generality along this line of investigations, especially when combined with some of the deeper results. These first slope results, however, are not discussed here since we are only interested in the much harder full Newton polygon in this paper. Similarly, some of the results in this paper could be improved somewhat in some cases if one restricts to the first few sides of the full Newton polygon.

We conclude this introductory chapter with some remarks on further questions. In this paper, our main concern is to determine when $f$ is ordinary, i.e., when the Newton polygon coincides with its lower bound. We have a reasonably good understanding of this question for a generic $f$ if $p$ is in the residue class $p \equiv 1 \pmod{D^*}$. What about other residue classes of $p$? In [p291, W3], there is a conjecture which classifies all large primes $p$ for which a generic $f$ is ordinary. This conjecture is true for any diagonal family, but is open in general. In connection with Grothendieck's specialization theorem, another general question, as asked by Mazur [Ma3] in certain geometric situations, is to study the geometry of the various stratifications of the parameter space imposed by the Newton polygon when one considers a family of polynomials over $\mathbf{F}_q$. The geometric case of polarized abelian varieties is understood by the work of Oort [Or]. The geometric case of curves (even plane curves) is already quite mysterious. Another interesting geometric case might be a family of Calabi-Yau varieties arising from mirror symmetry [Ya], generalizing the case of elliptic curves and the case of K3-surfaces. In our current exponential sum situation, the stratification question is already quite complicated even for the simplest family of one variable polynomials of degree $d$ (if $d$ is not too small). A related general question, as asked by Katz [Ka3] in the case of generalized Kloosterman polynomials, is to classify all the possible Newton polygons of a given class of non-diagonal polynomials $f$ when $p$ varies. For one variable polynomial $f(x) = x^3 + ax$ of degree 3, this was done by Sperber [S3] and there are already infinitely many (fortunately, well behaved) possibilities as $p$ varies. For one variable polynomials $f(x) = x^4 + ax^2 + bx$ of degree 4, see Hong [Ho] for some results. Further examples of this type, even for one variable polynomials of low degrees, would be of interest as they might suggest a possible general pattern. Still another question is to consider a fixed polynomial $f(x)$ with integer coefficients and investigate the possible density of ordinary primes for $f(x)$ in the set of all primes. In some cases, this density exists and is a positive rational number, see Corollary 4.3 for such a class of non-diagonal polynomials. In general, it is wide open. In the geometric case of an elliptic curve, Lang-Trotter [LT] has a precise conjecture on the distribution of supersingular primes, see Serre [Se] and Elkies [El] for some results on this conjecture. One can also consider the Newton polygon problem for a more general L-function such as Dwork's unit root L-function attached to a family of algebraic varieties over $\mathbf{F}_q$ or more generally a family of exponential sums. Such an L-function is deeper and no-longer rational in general. But it is $p$-adic meromorphic as Dwork [Dw2] conjectured, see [W7-8] for the proof of this conjecture. Very little is known about the Newton polygon for Dwork's unit root L-function. In the simplest elliptic family case, the Newton polygon problem is closely related to the Gouvêa-Mazur conjecture [GM] and the $p$-adic Ramanujuan-Peterson conjecture [W6] for modular forms, see Coleman [Co] and [W5] for some initial results in this direction.

## Chapter 2. A Lower Bound for the Newton Polygon

In this chapter, we review the lower bound of Adolphson-Sperber for the L-function of exponential sums. This will give us a suitable setting to describe precisely our basic $p$-adic question including the Adolphson-Sperber conjecture.

### 2.1. Non-degenerate Laurent polynomials.

Since $f$ is a Laurent polynomial, we can write $f$ as a finite sum of monomials:

$$f = \sum_{j=1}^{J} a_j x^{V_j}, \ a_j \neq 0,$$

where each $V_j = (v_{1j}, \cdots, v_{nj})$ is a lattice point in $\mathbf{Z}^n$ and the power $x^{V_j}$ simply means the product $x_1^{v_{1j}} \cdots x_n^{v_{nj}}$. Let $\Delta(f)$ be the convex closure in $\mathbf{R}^n$ generated by the origin and the lattice points $V_j$ ($1 \leq j \leq J$). This is called the Newton polyhedron of $f$. If $\delta$ is a subset of $\Delta(f)$, we define the restriction of $f$ to $\delta$ to be the Laurent polynomial

$$f^{\delta} = \sum_{V_j \in \delta} a_j x^{V_j}.$$

**Definition 2.1.** *The Laurent polynomial $f$ is called non-degenerate if for each (closed) face $\delta$ of $\Delta(f)$ of arbitrary dimension which does not contain the origin, the $n$ partial derivatives*

$$\{\frac{\partial f^{\delta}}{\partial x_1}, \cdots, \frac{\partial f^{\delta}}{\partial x_n}\}$$

*have no common zeros with $x_1 \cdots x_n \neq 0$ over the algebraic closure of $\mathbf{F}_q$.*

Note that the non-degenerate definition of $f$ depends only on the monomials of $f$ whose exponents are on the faces of $\Delta$ not containing the origin. Thus, if $g$ is another Laurent polynomial whose Newton polyhedron $\Delta(g)$ is contained in $\Delta(f)$ such that $\Delta(g)$ does not intersect the faces of $\Delta$ not containing the origin, then the sum $f + g$, which may be viewed as a deformation of $f$, is also non-degenerate whenever $f$ is non-degenerate.

For a simple example, let $f$ be the diagonal polynomial

$$f(x) = a_1 x_1^{d_1} + \cdots + a_n x_n^{d_n}, \ a_j \neq 0,$$

where the $d_j$ are positive integers. Then, the Newton polyhedron $\Delta(f)$ is the $n$-dimensional simplex in $\mathbf{R}^n$ whose vertices are the origin and the lattice points $(0, \cdots, d_j, \cdots, 0)$ for $1 \leq j \leq n$. It has only one face of co-dimension 1 not containing the origin. It is easy to check that if none of the $d_j$ is divisible by $p$, then $f$ is non-degenerate. If $g(x)$ is another polynomial such that each of its exponent $(v_1, \cdots, v_n)$ satisfies

$$\frac{v_1}{d_1} + \cdots + \frac{v_n}{d_n} < 1, \ v_j \geq 0,$$

then the sum $f(x) + g(x)$, which is no longer a diagonal polynomial, is also non-degenerate.

For another example, let $f(x)$ be the $n$-dimensional Kloosterman polynomial

$$f(x) = a_1 x_1 + \cdots + a_n x_n + a_{n+1} \frac{1}{x_1 \cdots x_n}, \ a_j \neq 0.$$

The Newton polyhedron $\Delta(f)$ is the $n$-dimensional simplex in $\mathbf{R}^n$ whose vertices are $(-1, \cdots, -1)$ and the lattice points $(0, \cdots, 1, \cdots, 0)$ for $1 \leq j \leq n$. It has $n+1$ faces of co-dimension 1, non of them containing the origin because the origin is an interior point of $\Delta(f)$. It is easy to check that $f(x)$ is non-degenerate for every prime number $p$. The following general result is due to Adolphson and Sperber [AS3].

**Theorem 2.1.** *Assume that $f$ is non-degenerate and $\Delta(f)$ is $n$-dimensional. Then, $\mathrm{L}^*(f,t)^{(-1)^{n-1}}$ is a polynomial of degree $n!\mathbf{V}(f)$, where $\mathbf{V}(f)$ is the volume of $\Delta(f)$. That is, there are $n!\mathbf{V}(f)$ non-zero algebraic integers $\alpha_i$ $(1 \leq i \leq n!\mathbf{V}(f))$ such that*

$$\mathrm{L}^*(f,t)^{(-1)^{n-1}} = \prod_{i=1}^{n!\mathbf{V}(f)} (1 - \alpha_i t). \tag{2.1}$$

With the same assumption on $f$, the complex absolute values of the reciprocal zeros $\alpha_i$ in (2.1) have been determined by Denef and Loeser [DL] using intersection cohomology and Deligne's fundamental results [De2] on Riemann hypothesis over finite fields, see also [De1] for earlier such applications to exponential sums. For each prime $\ell \neq p$, the reciprocal zeros $\alpha_i$ are always $\ell$-adic units by Deligne's integrality theorem. Our basic question is thus to determine the $p$-adic absolute values of the $n!\mathbf{V}(f)$ reciprocal zeros $\alpha_i$. Namely, we want to know

$$|\alpha_i|_p = ?? \quad \text{or} \quad \mathrm{ord}_q(\alpha_i) = ?? \tag{2.2}$$

for all $1 \leq i \leq n!\mathbf{V}(f)$. As indicated in the introduction, this is in general a very complicated problem. We have a complete understanding only in some special cases. An easy result is that $\mathrm{L}^*(f,t)^{(-1)^{n-1}}$ has exactly one reciprocal root $\alpha$ which is a $p$-adic unit (or slope zero). We are, however, interested in the harder question of determining the slopes of all reciprocal roots of the L-function. Before we move on, we need to reformulate the $n!\mathbf{V}(f)$ questions in (2.2) and pack them into a single question using the notion of Newton polygon to be defined in next section.

## 2.2. Definition of the Newton polygon.

Let $\Omega$ be the completion of an algebraic closure of the $p$-adic rational numbers $\mathbf{Q}_p$. Let $\mathbf{K}$ be a finite extension of $\mathbf{Q}_p$ contained in $\Omega$. Let $\pi$ be a fixed uniformizer of $\mathbf{K}$. Let

$$E(t) = 1 + \sum_{k \geq 1} a_k t^k, \ a_k \in \mathbf{K}, \ a_0 = 1$$

be a formal power series over $\mathbf{K}$ with constant term 1. Recall that $E(t)$ is called a $p$-adic entire function if $E(t)$ converges for all $t \in \Omega$. This is equivalent to the condition

$$\lim_{k \to \infty} \inf \frac{\operatorname{ord}_\pi(a_k)}{k} = +\infty. \tag{2.3}$$

In this case, the $p$-adic Weierstrass factorization theorem shows that there is a sequence $\alpha_k$ of elements in $\Omega$ such that

$$E(t) = \prod_{k \geq 1}(1 - \alpha_k t), \ \lim_{k \to \infty} \alpha_k = 0. \tag{2.4}$$

The elements $\alpha_k$ are the reciprocal zeros of $E(t)$. The factorization in (2.4) is the extension of the standard factorization theorem from a polynomial to a $p$-adic entire function. Note that unlike the complex case, there is no exponential factor in (2.4) because the exponential function is not a $p$-adic entire function.

**Definition 2.2.** *The Newton polygon $NP(E)$ of $E(t)$ is the convex closure from below in the plane $\mathbf{R}^2$ of the lattice points*

$$(0,0), (1, \operatorname{ord}_\pi a_1), \cdots, (k, \operatorname{ord}_\pi a_k), \cdots.$$

Assume that $E(t)$ is a $p$-adic entire function. Condition (2.3) implies that the slopes of the Newton polygon $NP(E)$ go to infinity. The $p$-adic absolute values (or the slopes) of the reciprocal zeros $\alpha_k$ can be easily read off from the Newton polygon $NP(E)$. In fact, we have the following well known fact which can be proved directly from the factorization in (2.4).

**Proposition 2.3.** *Let $E(t)$ be a $p$-adic entire function with constant term 1. The Newton polygon $NP(E)$ has a side of slope $s$ and horizontal length $h$ if and only if $E(t)$ has exactly $h$ reciprocal zeros $\alpha$ with slope $\operatorname{ord}_\pi \alpha = s$.*

Thus, to understand the $p$-adic absolute values of the zeros of a $p$-adic entire function $E(t)$, it suffices to determine the Newton polygon of $E(t)$. For the application to our basic problem of Section 2.1, our $p$-adic entire function will be the polynomial

$$E(t) = \mathrm{L}^*(f, t)^{(-1)^{n-1}}$$

of degree $n! \mathbf{V}(f)$ if $f$ is non-degenerate. We shall sometimes need to work directly on the chain level and thus work with entire functions instead of polynomials. This is in fact the approach taken in [W2-3], which has some flexibility to handle singular case as well. In practice, it is often possible to get a good lower bound for the Newton polygon $NP(E)$. But the precise determination of the Newton polygon is in general very difficult.

Note that instead of using the uniformizer $\pi$, sometimes it is more convenient to use a power of $\pi$ such as $p$ or $q$ in the definition of the Newton polygon. This does not change the substance and amounts to a rescaling of the Newton polygon.

### 2.3. A lower bound: the Hodge polygon.

We now return to the situation of L-functions as discussed in Section 2.1. We shall first describe the lower bound of Adolphson-Sperber for the Newton polygon of the polynomial $L^*(f,t)^{(-1)^{n-1}}$ for non-degenerate $f$. This lower bound depends in a simple way on the integral polyhedron $\Delta(f)$. It is called the Hodge polygon of $\Delta$ (or $f$) in analogue to the geometric case.

Let $\Delta$ denote the $n$-dimensional integral polyhedron $\Delta(f)$ in $\mathbf{R}^n$. Let $C(\Delta)$ be the cone in $\mathbf{R}^n$ generated by $\Delta$. Then $C(\Delta)$ is the union of all rays emanating from the origin and passing through $\Delta$. If $c$ is a real number, we define $c\Delta = \{cx|x \in \Delta\}$.

**Definition 2.4.** *For a point $u \in \mathbf{R}^n$, the weight $w(u)$ is defined to be the smallest non-negative real number $c$ such that $u \in c\Delta$. If such $c$ does not exist, we define $w(u) = \infty$.*

It is clear that $w(u)$ is finite if and only if $u \in C(\Delta)$. If $u \in C(\Delta)$ is not the origin, the ray emanating from the origin and passing through $u$ intersects $\Delta$ in a face $\delta$ of co-dimension 1 that does not contain the origin. The choice of the desired co-dimension 1 face $\delta$ is in general not unique unless the intersection point is in the interior of $\delta$. Let $\sum_{i=1}^{n} e_i X_i = 1$ be the equation of the hyperplane $\delta$ in $\mathbf{R}^n$, where the coefficients $e_i$ are uniquely determined rational numbers not all zero. Then, by standard arguments in linear programming, one finds that the weight function $w(u)$ can be computed using the formula:

$$w(u) = \sum_{i=1}^{n} e_i u_i, \tag{2.5}$$

where $(u_1, \cdots, u_n) = u$ denotes the coordinates of $u$. Let $D(\delta)$ be the least common denominator of the rational numbers $e_i$ $(1 \leq i \leq n)$. It follows from (2.5) that for a lattice point $u$ in $C(\delta)$, we have

$$w(u) \in \frac{1}{D(\delta)} \mathbf{Z}_{\geq 0}. \tag{2.6}$$

We leave it as an exercise to the reader to show that there are lattice points $u \in C(\delta)$ such that the denominator of $w(u)$ is exactly $D(\delta)$. That is, the denominator in (2.6) is optimal. Let $D(\Delta)$ be the least common multiple of all the $D(\delta)$, where $\delta$ runs over all the co-dimensional 1 faces of $\Delta$ which do not contain the origin. Then by (2.6), we deduce

$$w(\mathbf{Z}^n) \subseteq \frac{1}{D(\Delta)} \mathbf{Z}_{\geq 0} \cup \{+\infty\}, \tag{2.7}$$

where $\mathbf{Z}_{\geq 0}$ denotes the set of nonnegative integers. The integer $D(\Delta)$ is called the denominator of $\Delta$. It is the smallest positive integer for which (2.7) holds.

From the above explicit description of $w(u)$, one deduces the following two simple properties for the weight function $w(u)$.

**Lemma 2.5.** *For any real number $c \geq 0$, we have $w(cu) = cw(u)$.*

**Lemma 2.6.** *We have*
$$w(u + u') \leq w(u) + w(u'),$$
*with equality holding if and only if $u$ and $u'$ are co-facial, i.e., $w(u)^{-1}u$ and $w(u')^{-1}u'$ (assuming both $u$ and $u'$ are not the origin) lie on the same closed face of co-dimension 1 of $\Delta$ not containing the origin.*

For a non-negative integer $k$, let
$$W(k) = \text{card}\{u \in \mathbf{Z}^n \mid w(u) = \frac{k}{D(\Delta)}\}. \tag{2.8}$$
be the number of lattice points in $\mathbf{Z}^n$ with weight $k/D(\Delta)$. This is a finite number for each $k$. The integers $W(k)$ in (2.8) can be used to define a chain level Hodge polygon which is a lower bound [AS2] for the Newton polygon of the Fredholm determinant associated via Dwork's trace formula to the L-function. To define the Hodge polygon for the L-function itself, which is at cohomology level and which is our main concern here, we let
$$H(k) = \sum_{i=0}^{n}(-1)^i\binom{n}{i}W(k - iD(\Delta)). \tag{2.9}$$
This number is the number of lattice points of weight $k/D(\Delta)$ in a certain fundamental domain of $\Delta$ corresponding to a basis of the $p$-adic cohomology space used by Adolphson-Sperber to compute the L-function. Thus, $H(k)$ is a non-negative integer for each $k \in \mathbf{Z}_{\geq 0}$. Furthermore,
$$H(k) = 0, \text{ for } k > nD(\Delta)$$
and
$$\sum_{k=0}^{nD(\Delta)} H(k) = n!\mathbf{V}(\Delta), \tag{2.10}$$
where $\mathbf{V}(\Delta)$ is the volume of $\Delta$.

**Definition 2.7.** *The Hodge polygon $HP(\Delta)$ of $\Delta$ is defined to be the convex polygon in $\mathbf{R}^2$ with vertices $(0, 0)$ and*
$$\left(\sum_{k=0}^{m} H(k), \frac{1}{D(\Delta)}\sum_{k=0}^{m} kH(k)\right), \quad m = 0, 1, 2, \cdots, nD(\Delta). \tag{2.11}$$

The numbers $H(k)$ coincide with the Hodge numbers in the toric variety case in which case $D(\Delta) = 1$, see [AS4] for general results in this direction. This explains the term "Hodge polygon". Note that in geometric situation, one has $D(\Delta) = 1$ because one really considers the exponential sum of the new polynomial $x_{n+1}f(x)$ in order to study the zeta function of the hypersurface defined by $f$. The lower bound of Adolphson and Sperber [AS3] is then

**Theorem 2.8.** *Assume that $f$ is non-degenerate and $\Delta(f)$ is $n$-dimensional. Then, the Newton polygon of the polynomial $L^*(f,t)^{(-1)^{n-1}}$, computed with respect to $q$, lies above the Hodge polygon $HP(\Delta)$. Furthermore, the endpoints of the two polygons coincide.*

This theorem is the exponential sum analogue of the earlier results of Dwork-Mazur-Ogus for the zeta function of a smooth projective hypersurface [Dw1], a smooth projective liftable variety [Ma2] and a smooth proper variety with certain coefficients [BO] [Og2] in the direction of Katz's conjecture which says that the Newton polygon lies above the Hodge polygon, see also [AS4] for the geometric case of a smooth toric complete intersection. The Hodge polygon $HP(\Delta)$ is a geometric invariant of $f$. It depends only on the combinatorial property of the shape $\Delta(f)$ of the Laurent polynomial $f$. On the other hand, the Newton polygon of $L^*(f,t)^{(-1)^{n-1}}$ is a subtle arithmetic invariant of $f$, which depends not just on the combinatorial property of $\Delta(f)$, but also on subtle arithmetic property of $\Delta(f)$ and on the actual coefficients of the polynomial $f$. Thus, one cannot expect that there will be a simple general answer for the precise Newton polygon, although Grothendieck's specialization theorem gives a general structure theorem for the variation of the Newton polygon in an algebraic family.

## 2.4. The Adolphson-Sperber conjecture.

In the situation of Theorem 2.8, if the Newton polygon of the polynomial $L^*(f,t)^{(-1)^{n-1}}$ computed with respect to $q$ coincides with the Hodge polygon $HP(\Delta)$, then we say that $f$ is **ordinary**. In this case, the polynomial $L^*(f,t)^{(-1)^{n-1}}$ has exactly $H(k)$ reciprocal zeros $\alpha$ with $\mathrm{ord}_q(\alpha) = k/D(\Delta)$ for each $0 \leq k \leq nD(\Delta)$. If $f$ is not ordinary, then we say that $f$ is **non-ordinary**. The Hodge polygon is a reasonably good lower bound for the Newton polygon because the two polygons actually coincide in some cases. For example, if $f(x)$ is the one variable monomial $x^d$, then we have $D(\Delta) = d$ and $f(x)$ is ordinary if and only if $p \equiv 1 \pmod{D(\Delta)}$. For another example, let

$$f(x) = x_1 + \cdots + x_n + \frac{1}{x_1 \cdots x_n}$$

be the $n$-dimensional Kloosterman polynomial, then one checks that $D(\Delta) = 1$. Sperber's theorem says that $f(x)$ is ordinary for every prime $p$, i.e., for every prime $p \equiv 1 \pmod{D(\Delta)}$. These and many further examples treated by Sperber [S1-3] and Carpentier [Ca] led to the following general conjecture of Adolphson-Sperber, which suggests a closer relationship between the deeper arithmetic Newton polygon and the easier geometric (combinatorial) Hodge polygon.

**Conjecture 2.9.** *Let $\Delta$ be an $n$-dimensional integral convex polyhedron in $\mathbf{R}^n$ containing the origin. If*

$$p \equiv 1 \pmod{D(\Delta)}, \qquad (2.12)$$

*then the Newton polygon of $L^*(f,t)^{(-1)^{n-1}}$ coincides with the Hodge polygon $HP(\Delta)$ for a generic $f$ with $\Delta(f) = \Delta$, i.e., for all $f$ in a Zariski dense open subset of the parameter space $\mathcal{M}_p(\Delta)$ of non-degenerate $f$ with $\Delta(f) = \Delta$.*

The parameter space $\mathcal{M}_p(\Delta)$ is a variety defined over the prime field $\mathbf{F}_p$. It is the reduction modulo $p$ of a scheme $\mathcal{M}(\Delta)$ defined over $\mathbf{Z}$. For a point $f \in \mathcal{M}_p(\Delta)$ rational over some finite extension field $\mathbf{F}_q$, the Newton polygon of the L-function $L^*(f/\mathbf{F}_q,t)^{(-1)^{n-1}}$ computed with respect to $q$ is independent of the choice of the finite field $\mathbf{F}_q$ for which $f$ is defined. In particular, the ordinary property of $f$ is independent of the choice of the field $\mathbf{F}_q$ for which $f$ is defined. We can thus talk about the ordinary property of $f$ without specifying the field of definition for $f$. Only the prime number $p$ needs to be mentioned.

A simple ramification argument shows that in general the denominator $D(\Delta)$ in (2.12) cannot be replaced by any smaller positive integer. Thus, if Conjecture 2.9 is true, the condition in (2.12) is actually optimal in some sense. Our main question of this paper is to determine when the two polygons coincide, i.e., when the Laurent polynomial $f$ is ordinary. We have

**Theorem 2.10.** *Conjecture 2.9 is true in every low dimension $n \leq 3$ but false in every high dimension $n \geq 4$.*

This theorem shows that the full form of Conjecture 2.9 is in general false in every dimension $n \geq 4$. Counter-examples will be given in Section 3.4. Despite these negative results, we would like to know how far the conjecture is from being true and when it is true. It turns out that a weaker form of Conjecture 2.9 is always true. Namely, we have the following theorem from [W3].

**Theorem 2.11.** *Let $\Delta$ be an $n$-dimensional integral convex polyhedron in $\mathbf{R}^n$ containing the origin. There is an effectively computable finite positive integer $D^*(\Delta) \geq D(\Delta)$ such that if*

$$p \equiv 1 \ (\mathrm{mod}\ D^*(\Delta)), \tag{2.13}$$

*then the Newton polygon of $L^*(f,t)^{(-1)^{n-1}}$ coincides with the Hodge polygon $HP(\Delta)$ for a generic $f$ with $\Delta(f) = \Delta$, i.e., for all $f$ in a Zariski dense open subset of the parameter space $\mathcal{M}_p(\Delta)$ of non-degenerate $f$ with $\Delta(f) = \Delta$.*

In general, the smallest possible value for $D^*(\Delta)$ is quite subtle. It depends not just on the combinatorial shape of $\Delta$ but also on subtle arithmetic property of $\Delta$ as well. Again, there is probably no simple formula for the smallest $D^*(\Delta)$. Our method only gives a good effective upper bound in general which already depends on subtle arithmetic properties of $\Delta$. In many interesting higher dimensional cases, the full form of Conjecture 2.9 is also true, i.e., one can take $D^*(\Delta)$ to be $D(\Delta)$ in those cases.

The weaker form in Theorem 2.11 is trivial (if we do not care about the actual size of $D^*$) in the case that $\Delta$ is a simplex, because of the easy diagonal example and Grothendieck's specialization theorem. Thus, to prove Theorem 2.11 for a more general $\Delta$, a natural idea is to try to establish decomposition theorems for the generic Newton polygon corresponding to a suitable decomposition of $\Delta$. There is the facial decomposition theorem for the Newton polygon, first proved in [W2]. This decomposition simplifies $\Delta$ somewhat but is not enough to decompose $\Delta$ all the way into simplexes. In [W3], we introduced a degree polygon which is finer than the generic Newton polygon. We established two decomposition theorems for the degree polygon. One is called the star decomposition. The other is called the hyperplane decomposition. Using the method in [W3], we have recently proved a more flexible collapsing decomposition theorem for the degree polygon. The detail will be given elsewhere. Either the star decomposition or the better collapsing decomposition can be used to reduce the situation to the diagonal case. It should be noted that even for a simplex, our decomposition theorems (especially the new collapsing decomposition) can always be used to yield good upper bounds for the smallest possible $D^*(\Delta)$, which is far better (sometimes optimal) than what the trivial diagonal example and Stickelberger theorem would yield. This explains why our method can be used to prove the full form of Conjecture 2.9 in many cases including lower dimensional cases $n \leq 3$.

Mazur's conjecture for the zeta function of a generic hypersurface is another important example in higher dimensions. There should have many other similar geometric examples in the setting of toric hypersurfaces and toric complete intersections [AS4], which would follow from our hyperplane decomposition theorem and or the new collapsing decomposition theorem.

## Chapter 3. Diagonal Local Theory

In this chapter, we give a thorough treatment of the elementary diagonal case. Recall that a Laurent polynomial $f$ is called **diagonal** if $f$ has exactly $n$ non-zero terms and $\Delta(f)$ is $n$-dimensional (necessarily a simplex). In this case, the L-function can be computed explicitly by elementary method using Gauss sums, just as the classical case of the zeta function of a Fermat hypersurface. Thus, theoretically speaking or from algorithmic point of view, the Newton polygon in diagonal case can be completely determined using the Stickelberger theorem. It should, however, be noted that interesting arithmetic and combinatorial problems will often arise in the actual explicit calculation of the Newton polygon due to the diversity of the simplex $\Delta$. The diagonal case forms the building blocks for the general decomposition theory in Chapter 4. It also provides a rich source of explicit examples, including counter-examples of Conjecture 2.9.

### 3.1. $p$-action, Gauss sums and L-functions.

To describe the L-function in terms of Gauss sums, we first review the $p$-action on lattice points. Let $0, V_1, \cdots, V_n$ be the vertices of an $n$-dimensional integral simplex $\Delta$ in $\mathbf{R}^n$. Let $f$ be the diagonal Laurent polynomial

$$f(x) = \sum_{j=1}^{n} a_i x^{V_j}, a_j \in \mathbf{F}_q^*. \tag{3.1}$$

Let $M$ be the non-singular $n \times n$ matrix

$$M = (V_1, \cdots, V_n), \tag{3.2}$$

where each $V_j$ is written as a column vector. It is easy to check that $f$ is non-degenerate if $p$ is relatively prime to $\det M$. We assume that $p$ satisfies this condition. We consider the solutions of the following linear system

$$M \begin{pmatrix} r_1 \\ \vdots \\ r_n \end{pmatrix} \equiv 0 \ (\mathrm{mod}\, 1), \ r_i \text{ rational}, \ 0 \leq r_i < 1. \tag{3.3}$$

The map $(r_1, \cdots, r_n) \to r_1 V_1 + \cdots + r_n V_n$ clearly establishes a one-to-one correspondence between the solutions of (3.3) and the integral lattice points of the fundamental domain

$$\mathbf{R} V_1 + \cdots + \mathbf{R} V_n (\mathrm{mod}\ \mathbf{Z} V_1 + \cdots + \mathbf{Z} V_n). \tag{3.4}$$

Under this bijection, we can identify the solution set of (3.3) and the set of integral lattice points in the above fundamental domain. Let $S(\Delta)$ be the set of solutions $r$ of (3.3), which may be identified with the lattice points $u = Mr$ in the fundamental domain (3.4). It has a natural abelian group structure under addition modulo 1. By the theory of elementary abelian groups, the order of $S(\Delta)$ is precisely given by

$$|\det M| = n!\mathbf{V}(\Delta).$$

If $m$ is an integer relatively prime to the order of $S(\Delta)$, then multiplication by $m$ induces an automorphism of the finite abelian group $S(\Delta)$. This map is called the $m$-map (or $m$-action) of $S(\Delta)$, denoted by the notation $r \to \{mr\}$, where

$$\{mr\} = (\{mr_1\}, \cdots, \{mr_n\})$$

and $\{mr_i\}$ denotes the fractional part of the real number $mr_i$. For each element $r \in S(\Delta)$, let $d(m,r)$ be the smallest positive integer such that multiplication by $m^{d(m,r)}$ acts trivially on $r$, i.e.,

$$(m^{d(m,r)} - 1)r \in \mathbf{Z}^n.$$

The integer $d(m,r)$ is the order of the $m$-map restricted to the cyclic subgroup generated by $r$. It is called the $m$-degree of $r$. For each positive integer $d$, let $S(m,d)$ be the set of $r \in S(\Delta)$ such that $d(m,r) = d$. Of course, the set $S(m,d)$ is empty for all large $d$ since $S(\Delta)$ is finite. We have the disjoint $m$-degree decomposition

$$S(\Delta) = \cup_{d \geq 1} S(m,d).$$

Since $p$ is relatively prime to the order of $S(\Delta)$, we will use the case that $m$ is a power of $p$, such as $p$ and $q$.

In order to compute the L-function, we now recall the definition of Gauss sums. Let $\chi$ be the Teichmüller character of the multplicative group $\mathbf{F}_q^*$. For $a \in \mathbf{F}_q^*$, the value $\chi(a)$ is just the $(q-1)$-th root of unity in the $p$-adic field $\Omega$ such that $\chi(a)$ modulo $p$ reduces to $a$. Define the $(q-1)$ Gauss sums over $\mathbf{F}_q$ by

$$G_k(q) = -\sum_{a \in \mathbf{F}_q^*} \chi(a)^{-k} \zeta_p^{\mathrm{Tr}(a)} \quad (0 \leq k \leq q-2),$$

where Tr denotes the trace map from $\mathbf{F}_q$ to the prime field $\mathbf{F}_p$. For each $a \in \mathbf{F}_q^*$, the Gauss sums satisfy the following interpolation relation

$$\zeta_p^{\mathrm{Tr}(a)} = \sum_{k=0}^{q-2} \frac{G_k(q)}{1-q} \chi(a)^k.$$

One then calculates that

$$\begin{aligned}
S_1^*(f) &= \sum_{x_j \in \mathbf{F}_q^*} \zeta_p^{\mathrm{Tr}(f(x))} \\
&= \sum_{x_j \in \mathbf{F}_q^*} \prod_{i=1}^n \zeta_p^{\mathrm{Tr}(a_i x^{V_i})} \\
&= \sum_{x_j \in \mathbf{F}_q^*} \prod_{i=1}^n \sum_{k_i=0}^{q-2} \frac{G_{k_i}(q)}{1-q} \chi(a_i)^{k_i} \chi(x^{V_i})^{k_i} \qquad (3.5)\\
&= \sum_{k_1=0}^{q-2} \cdots \sum_{k_n=0}^{q-2} (\prod_{i=1}^n \frac{G_{k_i}(q)}{1-q} \chi(a_i)^{k_i}) \sum_{x_j \in \mathbf{F}_q^*} \chi(x^{k_1 V_1 + \cdots + k_n V_n}) \\
&= (-1)^n \sum_{k_1 V_1 + \cdots + k_n V_n \equiv 0 (\bmod q-1)} \prod_{i=1}^n \chi(a_i)^{k_i} G_{k_i}(q).
\end{aligned}$$

This gives a formula for the exponential sum $S_1^*(f)$ over the field $\mathbf{F}_q$. Replacing $q$ by $q^k$, one gets a formula for the exponential sum $S_k^*(f)$ over the $k$-th extension field $\mathbf{F}_{q^k}$.

If $r = (r_1, \cdots, r_n)$ and $r' = (r'_1, \cdots, r'_n)$ are two elements in $S(p, d)$ which are in the same orbit under the $p$-action, that is

$$r' = \{p^k r\}$$

for some integer $k$, then one checks from the above definition of Gauss sums that

$$G_{r_i(q^d-1)}(q^d) = G_{r'_i(q^d-1)}(q^d)$$

for all $1 \leq i \leq n$, where $G_k(q^d)$ is the Gauss sum defined over the finite extension field $\mathbf{F}_{q^d}$. Thus, the $p$-action and the $q$-action do not change the Gauss sum. For $r \in S(q, d)$, the well known Hasse-Davenport relation says that

$$G_{r(q^{dk}-1)}(q^{dk}) = G_{r(q^d-1)}(q^d)^k,$$

for every positive integer $k$. Using this relation and (3.5), one then obtains the following explicit formula [W2] for the L-function.

**Theorem 3.1.** *Let $p$ be relatively prime to $\det M$. Using the above notations, we have*

$$L^*(f/\mathbf{F}_q, t)^{(-1)^{n-1}} = \prod_{d \geq 1} \prod_{r \in S(q,d)} \left(1 - t^d \prod_{i=1}^n \chi(a_i)^{r_i(q^d-1)} G_{r_i(q^d-1)}(q^d)\right)^{\frac{1}{d}}, \tag{3.6}$$

*where $r = (r_1, \cdots, r_n)$.*

Note that we used the $q$-action in the definition of the set $S(q, d)$ which occurs in (3.6). For each of the $d$ points in the orbit of $r \in S(q, d)$ under the $q$-action, the corresponding factor in (3.6) is the same. Thus, we can remove the power $1/d$ if we restrict $r$ to run over the $q$-orbits (closed points) of $S(q, d)$. In particular, this shows that the right side of (3.6) is indeed a polynomial of degree $n!\mathbf{V}(\Delta)$.

### 3.2. Applications of Stickelberger's theorem.

The Stickelberger theorem for Gauss sums can be described as follows.

**Theorem 3.2.** *Let $0 \leq k \leq q - 2$. Let $\sigma_p(k)$ be the sum of the $p$-digits of $k$ in its base $p$ expansion. That is,*

$$\sigma_p(k) = k_0 + k_1 + k_2 + \cdots, \quad k = k_0 + k_1 p + k_2 p^2 + \cdots, \quad 0 \leq k_i \leq p - 1.$$

*Then,*

$$\mathrm{ord}_p G_k(q) = \frac{\sigma_p(k)}{p-1}.$$

Combining Theorems 3.1 and 3.2, one can then completely determine the $p$-adic absolute values of the reciprocal zeros of $L^*(f, t)^{(-1)^{n-1}}$. In particular, the Newton polygon is independent of the non-zero coefficients $a_j$ of the diagonal Laurent polynomial $f$. Thus, for simplicity, we may assume that all the coefficients are 1. Namely,

$$f = \sum_{j=1}^n x^{V_j}. \tag{3.7}$$

Applying Theorem 3.1 to the polynomial in (3.7) and $q = p$, we obtain

**Corollary 3.3.** *Let $p$ be relatively prime to $\det M$. For $f$ in (3.7), we have*

$$L^*(f/\mathbf{F}_p, t)^{(-1)^{n-1}} = \prod_{d \geq 1} \prod_{r \in S(p,d)} \left(1 - t^d \prod_{i=1}^n G_{r_i(p^d-1)}(p^d)\right)^{\frac{1}{d}}, \tag{3.8}$$

where $r = (r_1, \cdots, r_n)$ and the $p$-action is used.

The Newton polygon computed with respect to $q$ of the polynomial in Theorem 3.1 is the same as the Newton polygon computed with respect to $p$ of the polynomial in Corollary 3.3. This follows directly from the relationship between the two L-functions: the reciprocal zeros of $L^*(f/\mathbf{F}_q, t)^{(-1)^{n-1}}$ for the polynomial $f$ in (3.7) are exactly the $k$-th power of the reciprocal zeros of $L^*(f/\mathbf{F}_p, t)^{(-1)^{n-1}}$, where $q = p^k$.

If now $r = (r_1, \cdots, r_n)$ is an element of $S(p,d)$ and $\alpha_r$ is any one of the $d$ reciprocal roots (which differ only by a $d$-th root of unity) of the corresponding factor in (3.8), then the Stickelberger theorem implies that

$$\begin{aligned}
\operatorname{ord}_p(\alpha_r) &= \frac{1}{d} \operatorname{ord}_p \prod_{i=1}^n G_{r_i(p^d-1)}(p^d) \\
&= \frac{1}{d(p-1)} \sum_{i=1}^n \sigma_p(r_i(p^d - 1)) \\
&= \frac{1}{d(p-1)} \sum_{i=1}^n \frac{p-1}{p^d-1} \sum_{j=0}^{d-1} \{p^j r_i\}(p^d - 1) \\
&= \frac{1}{d} \sum_{j=0}^{d-1} |\{p^j r\}|,
\end{aligned} \tag{3.9}$$

where for $r \in S(\Delta)$, the norm $|r| = r_1 + \cdots + r_n$ is just the weight $w(u)$ of the corresponding lattice point $u = Mr$ in the fundamental domain (3.4). From (3.9), one can write down a (complicated) formula for the Newton polygon of the diagonal L-function.

### 3.3. Ordinary criterion and ordinary primes.

To determine when the diagonal Laurent polynomial $f$ in (3.1) is ordinary, we need to describe the Hodge polygon in the diagonal case. For this purpose, we consider the following linear system for $u \in \mathbf{Z}^n$,

$$M \begin{pmatrix} r_1 \\ \vdots \\ r_n \end{pmatrix} = u, \ 0 \leq r_i. \tag{3.10}$$

If $u \in C(\Delta)$, system (3.10) has exactly one solution $r = (r_1, \cdots, r_n)$. In this case, the weight $w(u)$ of $u$ is given by

$$w(u) = r_1 + \cdots + r_n = |r|.$$

Thus, $W(k)$ is the number of lattice points $u \in \mathbf{Z}^n$ such that

$$w(u) = r_1 + \cdots + r_n = \frac{k}{D}. \tag{3.11}$$

By the inclusion-exclusion principle, the Hodge number
$$H(k) = W(k) - \binom{n}{1}W(k-D) + \binom{n}{2}W(k-2D) - \cdots$$
is simply the number of lattice points in the fundamental domain (3.4) with weight $k/D$. The Hodge polygon $HP(\Delta)$ is the polygon with vertices $(0,0)$ and
$$(\sum_{k=0}^{m} H(k), \sum_{k=0}^{m} \frac{k}{D}H(k)), \ m = 0, 1, \cdots, nD.$$

**Theorem 3.4.** *The diagonal Laurent polynomial $f$ in (3.1) is ordinary over $\mathbf{F}_q$ if and only if the norm function $|r|$ on $S(\Delta)$ is stable under the $p$-action: That is, for each $r \in S(\Delta)$, we have*
$$|r| = |\{pr\}|.$$
*Equivalently, this means that the weight function $w(u)$ on the lattice points of the fundamental domain (3.4) is stable under the $p$-action:*
$$w(u) = w(\{pu\}). \tag{3.12}$$

**Proof.** If the norm function $|r|$ on $S(\Delta)$ is stable under the $p$-action, then for any $r \in S(p,d)$, (3.9) shows that each corresponding reciprocal root $\alpha_r$ satisfies
$$\mathrm{ord}_p(\alpha_r) = |r| = w(u),$$
where $u = Mr$. It follows from Corollary 3.3 that the L-function $L^*(f/\mathbf{F}_p, t)^{(-1)^{n-1}}$ has exactly $H(k)$ reciprocal roots $\alpha$ with slope $\mathrm{ord}_p(\alpha) = k/D$. This means that $f$ is ordinary at $p$.

Conversely, if $f$ is ordinary at $p$, one shows by induction that each reciprocal root $\alpha_r$ corresponding to $r \in S(p,d)$ satisfies
$$\mathrm{ord}_p(\alpha_r) = |r|.$$
Since the left side of this equation depends only on the $p$-orbit of $r$, we must have
$$|r| = |\{pr\}| = |\{p^2 r\}| = \cdots.$$
This last equation is satisfied if and only if the norm function $|r|$ on $S(\Delta)$ is stable under the $p$-action. The theorem is proved.

In particular, if the $p$-action on $S(\Delta)$ is trivial, then the weight function $w(u)$ on the lattice points of the fundamental domain (3.4) is automatically stable under the $p$-action and hence the diagonal Laurent polynomial $f$ in (3.1) is ordinary over $\mathbf{F}_q$. Let $d_1 \mid d_2 \mid \cdots \mid d_n$ be the invariant factors of the matrix $M$ (or the finite abelian group $S(\Delta)$). That is, there are two matrices $P, Q \in \mathrm{GL}_n(\mathbf{Z})$ such that
$$PMQ = \mathrm{diag}(d_1, d_2, \cdots, d_n)$$
is a diagonal matrix with $d_i | d_{i+1}$. By the theory of finite abelian groups, multiplication by the largest invariant factor $d_n$ kills the finite abelian group $S(\Delta)$. Thus, if $p-1$ is divisible by $d_n$, the $p$-action on $S(\Delta)$ becomes trivial. We obtain

**Corollary 3.5.** *Let $d_n$ be the largest invariant factor of $M$. If $p \equiv 1 \pmod{d_n}$, then the diagonal Laurent polynomial in (3.1) is ordinary at $p$.*

As an application, we get information about the variation of the Newton polygon when $p$ varies. For example, we have

**Corollary 3.6.** *Let $f(x)$ be a diagonal Laurent polynomial with integer coefficients. Let $d_n$ be the largest invariant factor of the $n$-dimensional simplex $\Delta(f)$. For a positive real number $t$, let $\pi_f(t)$ be the number of primes $p \leq t$ such that $f(x)$ is ordinary at $p$. Then, there is a positive integer $\mu(\Delta) \leq \varphi(d_n)$ such that the following asymptotic formula holds*

$$\pi_f(t) \sim \frac{\mu(\Delta)}{\varphi(d_n)} \frac{t}{\log t},$$

*where $\varphi(d_n)$ denotes the Euler function.*

In fact, if $p_1$ is an ordinary prime for $f$, then Theorem 3.4 implies that $f$ is ordinary for every prime $p \equiv p_1 \pmod{d_n}$. Thus, the set of ordinary primes for $f$ are exactly the primes in certain residue classes modulo $d_n$. The integer $\mu(\Delta)$ is the number of positive integers $m \leq d_n$ such that $(m, d_n) = 1$ and the weight function $w(u)$ on the fundamental domain of $\Delta(f)$ is stable under the $m$-action.

We now discuss the relationship between the denominator $D(\Delta)$ and the largest invariant factor $d_n$. By definition, the denominator $D(\Delta)$ is the least common denominator for the coordinates of the vector $(e_1, \cdots, e_n)$ which is the unique solution of the linear system

$$(e_1, \cdots, e_n)M = (1, \cdots, 1).$$

Solving this linear system, one finds that

$$(e_1, \cdots, e_n) = (1, \cdots, 1)M^{-1} \in \frac{1}{d_n}\mathbf{Z}^n. \tag{3.13}$$

This shows that $D(\Delta)$ is a factor of $d_n$ and we have the inequality

$$D(\Delta) \leq d_n.$$

In general, $D(\Delta)$ will be a proper factor of $d_n$ if $n \geq 2$. The simplest example is to take the 2 variable diagonal Laurent polynomial $f(x) = x_1^d x_2^{1-d} + x_2$. One checks that $D = 1$ but the largest invariant factor $d_2 = d$.

We finish this section by showing that the full form of the Adolphson-Sperber conjecture holds in the indecomposable case for $n \leq 3$. Let $\delta$ be the $(n-1)$-dimensional face generated by the $V_j$ ($1 \leq j \leq n$). This is the unique co-dimension 1 face of the simplex $\Delta$ not containing the origin. We say that $\Delta$ is **indecomposable** if the face $\delta$ contains no lattice points other than the vertices $V_j$. Note that the larger $\Delta$ may have non-vertex lattice points with weight less than 1.

**Corollary 3.7.** *Conjecture 2.9 is true for indecomposable $\Delta$ if $n \leq 3$.*

**Proof.** By Corollary 3.5 and the deformation consequence of Theorem 4.1 in Section 4.1, we only need to prove that the denominator $D(\Delta)$ is equal to the largest invariant factor $d_n$ if $n \leq 3$. Let

$$e_1 X_1 + \cdots + e_n X_n = D, \quad \gcd(e_1, \cdots, e_n) = 1 \tag{3.14}$$

be the equation of the hyperplane $\delta$, where the $e_i$ are relatively prime integers. The Euclidean algorithm shows that there is a matrix $P \in \mathrm{GL}_n(\mathbf{Z})$ such that the first row of $P$ is $(e_1, \cdots, e_n)$. Since each column vector $V_j$ of $M$ satisfies equation (3.14), we see that the first row of the product matrix $PM$ is $(D, \cdots, D)$. This normalization shows that without loss of generality, we may assume that $\delta$ is defined by the equation $X_1 = D$. Recall that

$$D \leq d_n \leq |\det(M)|.$$

Thus, it suffices to prove that $\det(M) = \pm D$.

If $n = 1$, it is trivial that $D = d_1$. If $n = 2$, then $\delta$ is the line segment between $V_1 = (D, v_1)$ and $V_2 = (D, v_2)$ for some integers $v_1 \neq v_2$. Since there is no lattice point strictly between $v_1$ and $v_2$, we deduce that $|v_1 - v_2| = 1$ and thus $\det(M) = \pm D$.

If $n = 3$, let $\delta_1$ be the 2-dimensional integral simplex in $\mathbf{R}^2$ whose vertices are the origin, $U_2$ and $U_3$, where $U_i$ ($2 \leq i \leq 3$) is the last two coordinates of $V_i - V_1$. Since there are no non-vertex lattice points on $\delta$, there are also no non-vertex lattice points on $\delta_1$. But $\delta_1$ has dimension 2. We deduce the stronger property that there are no non-vertex lattice points in the fundamental domain $\mathbf{R}^2(\mathrm{mod}\mathbf{Z}U_2 + \mathbf{Z}U_3)$. In fact, if $r_2U_2 + r_2U_3$ ($0 < r_i < 1$) were a non-vertex lattice point in the fundamental domain of $\delta_1$ but not in $\delta_1$, then $r_2 + r_3 > 1$ and we deduce that its mirror $(1 - r_2)U_2 + (1 - r_3)U_3$ is an interior lattice point of $\delta_1$ since $(1 - r_2) + (1 - r_3) < 1$. No non-vertex lattice point in the fundamental domain of the simplex $\delta_1$ implies that $\det(U_2, U_3) = \pm 1$. It follows that

$$\det(V_1, V_2, V_3) = \det(V_1, V_2 - V_1, V_3 - V_1) = D\det(U_2, U_3) = \pm D.$$

The proof is complete.

### 3.4. Counter-examples in high dimensions

The discussion in the previous section shows that, in searching for diagonal counter-examples of Conjecture 2.9, we should look at those $\Delta$ with $D(\Delta) < d_n$. Furthermore, to insure that the diagonal family in (3.1) is already the generic family, we need to choose $\Delta$ to be **minimal**, i.e., there are no other lattice points on $\Delta$ other than the vertices. This is the approach taken in [W3]. The following 5-dimensional counter-example is given in [W3].

Let $\Delta$ be the 5-dimensional simplex in $\mathbf{R}^5$ whose non-zero vertices are the column vectors $V_i$ ($1 \leq i \leq 5$) of the following matrix:

$$M = \begin{pmatrix} 1 & 1 & 1 & 1 & 1 \\ 0 & 0 & 1 & 1 & 1 \\ 0 & 1 & 0 & 1 & 1 \\ 0 & 1 & 1 & 0 & 1 \\ 0 & 1 & 1 & 1 & 0 \end{pmatrix}$$

Clearly, $D(\Delta) = 1$ since the first coordinate of each $V_i$ is 1. One computes that the determinant of the matrix $M$ has absolute value 3. A direct proof shows that there are no lattice points on $\Delta$ except for the vertices. Thus, the generic family is the diagonal family. In particular, $\mathrm{L}^*(f, t)$ is a polynomial of degree 3 if $p \neq 3$. There are only three lattice points contained in the fundamental domain $\mathbf{R}^5/\mathbf{Z}V_1 + \cdots + \mathbf{Z}V_5$. One of them is the origin. The other two are

$$(2,1,1,1,1) = \frac{2}{3}V_1 + \sum_{i=2}^{5}\frac{1}{3}V_i, \quad (3,2,2,2,2) = \frac{1}{3}V_1 + \sum_{i=2}^{5}\frac{2}{3}V_i. \tag{3.15}$$

These two lattice points have weight 2 and 3, respectively. One checks that if $p \equiv 2 \pmod{3}$, the $p$-action permutes the two points in (3.15) and thus the weight function is not stable under the $p$-action. We conclude that the Newton polygon lies strictly above its lower bound if $p \equiv 2 \pmod{3}$.

To construct counter-examples in higher dimensions ($n \geq 6$), we can simply let $\Delta$ be the integral simplex in $\mathbf{R}^n$ generated by the origin and the column vectors of the $n \times n$ matrix

$$M^* = \begin{pmatrix} M & B \\ 0 & I_{n-5} \end{pmatrix}, \tag{3.16}$$

where $M$ is the above $5 \times 5$ matrix, $I_{n-5}$ is the identity matrix of order $(n-5)$, the first row of $B$ is $(1, 1, \cdots, 1)$ and all other rows of $B$ are zero. The first row of $M^*$ is $(1, 1, \cdots, 1)$ and so $D(\Delta) = 1$. One checks that there are no lattice points on $\Delta$ other than the vertices. The generic family is indeed the diagonal family. However, if $p \equiv 2 \pmod{3}$, the above argument shows that the Newton polygon of $L^*(f,t)^{(-1)^{n-1}}$ lies strictly above its lower bound.

Counter-examples in dimension 4 are a little harder to construct and to prove. Let $D \geq 2$ and $k \geq 2$ be any two given positive integers. Let $f$ be the 4-variable diagonal Laurent polynomial

$$f(x) = a_1 x_1^D + a_2 x_1^D x_2 + a_3 x_1^D x_2 x_3 + a_4 x_1^D x_3^{-1} x_4^{D^k}, \tag{3.17}$$

where the coefficients are non-zero elements of $\mathbf{F}_q$. Thus, $\Delta = \Delta(f)$ is the 4-dimensional simplex in $\mathbf{R}^4$ whose non-zero vertices are the column vectors of the following matrix

$$M = \begin{pmatrix} D & D & D & D \\ 0 & 1 & 1 & 0 \\ 0 & 0 & 1 & -1 \\ 0 & 0 & 0 & D^k \end{pmatrix}.$$

It is clear that $D(\Delta) = D$, $d_4 = D^k$ and $\det(M) = D^{k+1}$. The lattice points in the fundamental domain $\mathbf{R}^4/\mathbf{Z}V_1 + \cdots + \mathbf{Z}V_4$ are of the form

$$\frac{j_1}{D^k}V_1 + \{\frac{D^k - j_4}{D^k}\}V_2 + \frac{j_4}{D^k}V_3 + \frac{j_4}{D^k}V_4,$$

where $\{x\}$ denotes the fractional part of the real number $x$, $j_1$ and $j_4$ are integers such that

$$0 \leq j_1, j_4 \leq D^k - 1, \quad j_1 + j_4 \equiv 0 \pmod{D^{k-1}}. \tag{3.18}$$

There are exactly $D^{k+1}$ choices in (3.18), corresponding to the number of lattice points in the fundamental domain. There are exactly $D$ lattice points in the fundamental domain with weight smaller than 1 corresponding to the case $j_4 = 0$. There are no lattice points on $\Delta$ with weight 1 other than the vertices. The next smallest possible weight for a lattice point in the fundamental domain is $1 + \frac{1}{D}$. Take $j_4 = 1$ and $j_1 = D^{k-1} - 1$, then

$$u = (D+1, 1, 0, 1) = \frac{D^{k-1} - 1}{D^k} V_1 + \frac{D^k - 1}{D^k} V_2 + \frac{1}{D^k} V_3 + \frac{1}{D^k} V_4$$

is a lattice point in the fundamental domain with weight $1 + \frac{1}{D}$. The $p$-action of this point is given by

$$\{pu\} = \{\frac{p(D^{k-1} - 1)}{D^k}\} V_1 + \{\frac{p(D^k - 1)}{D^k}\} V_2 + \{\frac{p}{D^k}\} V_3 + \{\frac{p}{D^k}\} V_4. \tag{3.19}$$

The weight of this point is bounded from below by

$$w(\{pu\}) \geq 1 + \{\frac{p}{D^k}\} \tag{3.20}$$

since the sum of the middle two coefficients in (3.19) is already 1. This shows that the weight function is not stable under the $p$-action if

$$\{\frac{p}{D^k}\} > \frac{1}{D}.$$

The last condition is satisfied, for instance, if we take $p$ such that $p - (1 + D^{k-1})$ is divisible by $D^k$. Such $p$ satisfies the condition that $p - 1$ is divisible by $D^{k-1}$ (in particular by $D$) but not divisible by $D^k$. Thus, the diagonal Laurent polynomial $f$ in (3.17) is not ordinary for such $p$.

Although the Laurent polynomial in (3.17) is not yet generic in terms of $\Delta$. It is "essentially" generic as far as the ordinary property of $f$ is concerned. This is because there are no lattice points on $\Delta$ with weight 1 other than the non-zero vertices. Thus, the diagonal Laurent polynomial in (3.17) is the "leading form" of the generic Laurent polynomial with respect to $\Delta$. The only terms we missed are those terms with weight strictly less than 1. Theorem 4.1 in next section shows that deformations by such "error terms" (with weight less than 1) have no effect on the ordinary property of $f$, although it would indeed change the Newton polygon in non-ordinary case. This shows that Conjecture 2.9 is false in dimension 4 as well. This counter-example together with the trick in (3.16) gives the following result.

**Corollary 3.8.** *Let $D^*(\Delta)$ be the smallest positive integer such that Conjecture 2.9 is true for all $p$ in the residue class $p \equiv 1 \pmod{D^*(\Delta)}$. Then for each integer $n \geq 4$, we have*

$$\lim_{\dim(\Delta) = n} \sup \frac{D^*(\Delta)}{D(\Delta)} = +\infty.$$

This shows that for each integer $n \geq 4$, the two numbers $D^*(\Delta)$ and $D(\Delta)$ can differ as much as one would like. Nevertheless, the full form of Conjecture 2.9 is true in many important cases as we shall see.

## Chapter 4. Global Decomposition Theory

The diagonal case is in principal well understood by the discussion in the previous chapter. We now turn to the general case and describe several decomposition theorems for the Newton polygon of the L-function $L^*(f, t)^{(-1)^{n-1}}$. Except for the facial decomposition theorem, all other decomposition theorems are originally proved for the finer degree polygon which we have defined. To describe these purely in terms of generic Newton polygons, the results become weaker. However, they are sufficient for our applications here. The star decomposition theorem [W3] is not included since it is not used in our chosen applications here.

### 4.1. Facial decomposition for the Newton polygon.

The facial decomposition theorem described in this section allows us to determine when $f$ is ordinary in certain non-diagonal cases.

Let $f(x)$ be a Laurent polynomial over $\mathbf{F}_q$ such that $\Delta(f) = \Delta$ is $n$-dimensional. We assume that $f$ is non-degenerate. Let $\delta_1, \cdots, \delta_k$ be all the co-dimension 1 faces of $\Delta$ which do not contain the origin. Let $f^{\delta_i}$ be the restriction of $f$ to the face $\delta_i$. Then, $\Delta(f^{\delta_i}) = \Delta_i$ is $n$-dimensional. Furthermore, $f$ is non-degenerate if and only if $f^{\delta_i}$ is non-degenerate for all $1 \leq i \leq k$. Since we assumed that $f$ is non-degenerate, it follows that each $f^{\delta_i}$ is also non-degenerate.

The following facial decomposition theorem is taken from [W2].

**Theorem 4.1 (facial decomposition).** *Let $f$ be non-degenerate and let $\Delta(f)$ be $n$-dimensional. Then $f$ is ordinary over $\mathbf{F}_q$ if and only if each $f^{\delta_i}$ ($1 \leq i \leq k$) is ordinary over $\mathbf{F}_q$. Equivalently, $f$ is non-ordinary over $\mathbf{F}_q$ if and only if some $f^{\delta_i}$ is non-ordinary over $\mathbf{F}_q$.*

This theorem shows that as far as the ordinary property of $f$ is concerned, we may assume that $\Delta(f)$ has only one face of co-dimension 1 not containing the origin. Combining Theorem 4.1 with the results of Chapter 3 (Theorem 3.4), we obtain

**Corollary 4.2.** *If each $f^{\delta_i}$ ($1 \leq i \leq k$) is a diagonal Laurent polynomial, then we have a complete classification of the primes $p$ for which $f$ is ordinary.*

Together with Corollary 3.6, we deduce

**Corollary 4.3.** *Let $f(x)$ be a Laurent polynomial with integer coefficients and with $\Delta(f) = \Delta$. For a positive real number $t$, let $\pi_f(t)$ be the number of primes $p \leq t$ such that $f(x)$ is ordinary at $p$. Assume that $f^{\delta_i}$ is diagonal for each $1 \leq i \leq k$. Then, there is a positive rational number $r(\Delta) \leq 1$ such that we have the following asymptotic formula*
$$\pi_f(t) \sim r(\Delta) \frac{t}{\log t}.$$

Combining Theorem 4.1 with Corollary 3.5, we deduce

**Corollary 4.4.** *If each $f^{\delta_i}$ ($1 \leq i \leq k$) is a diagonal Laurent polynomial whose largest invariant factor is $d_{ni}$, then $f$ is ordinary at $p$ whenever $p - 1$ is divisible by the least common multiple $\mathrm{lcm}(d_{n1}, \cdots, d_{nk})$.*

For a simple example with $k = 1$, let $f(x)$ be a polynomial over $\mathbf{F}_q$ of the form

$$f(x) = a_1 x_1^d + a_2 x_2^d + \cdots + a_n x_n^d + g(x), \ a_j \in \mathbf{F}_q^*, \tag{4.1}$$

where $d$ is a positive integer and $g(x)$ is any polynomial of degree smaller than $d$. The polynomial $f(x)$ in (4.1) is a deformation of the Fermat polynomial. It is non-degenerate if and only if $d$ is not divisible by $p$, in which case the L-function $L^*(f, T)^{(-1)^{n-1}}$ is a polynomial of degree $d^n$. Theorem 4.1 shows that $f(x)$ is ordinary over $\mathbf{F}_q$ if and only if the leading diagonal form $f^\delta = a_1 x_1^d + \cdots + a_n x_n^d$ is ordinary over $\mathbf{F}_q$. This last condition holds if and only if $p - 1$ is divisible by $d$.

For a simple example with $k > 1$, let $f(x)$ be the $n$-dimensional generalized Kloosterman polynomial

$$f(x) = a_1 x_1 + \cdots + a_n x_n + a_{n+1} \frac{1}{x_1^{v_1} \cdots x_n^{v_n}}, \quad a_j \neq 0, \tag{4.2}$$

where each $v_j$ is a positive integer. The polynomial $f(x)$ in (4.2) is non-degenerate if and only if none of the $v_j$ is divisible by $p$, in which case the L-function $L^*(f, T)^{(-1)^{n-1}}$ is a polynomial of degree $1 + v_1 + \cdots + v_n$. Let $\delta_i$ be the $(n-1)$-dimensional simplex formed by all the exponents of $f(x)$ with the $i$-th exponent removed, where $1 \leq i \leq n+1$. Let $\Delta_i$ be the $n$-dimensional simplex generated by the origin and $\delta_i$. The invariant factors of $\Delta_{n+1}$ are all 1. Thus, $f^{\delta_{n+1}}$ is ordinary for all $p$. For $1 \leq i \leq n$, the invariant factors of $\Delta_i$ are given by

$$\{1, \cdots, 1, v_i\},$$

and thus $f^{\delta_i}$ is ordinary if $p - 1$ is divisible by $v_i$. We conclude that the generalized Kloosterman polynomial $f(x)$ in (4.2) is ordinary for all $p$ such that $p - 1$ is divisible by $\mathrm{lcm}(v_1, \cdots, v_n)$. This result was first proved by Sperber [S2], see also [S1] for the classical case $v_1 = \cdots v_n = 1$. We leave it as an exercise to the reader to apply the above method to study the more general polynomial

$$f(x) = a_1 x_1^{u_1} + \cdots + a_n x_n^{u_n} + a_{n+1} x_1^{-v_1} \cdots x_n^{-v_n}, \quad a_j \neq 0, u_i, v_j > 0, \tag{4.3}$$

which will include the results of Carpentier [Ca] as a special case. Our method actually gives more for the polynomial in (4.2)-(4.3), since Corollary 4.2 shows that we have a complete classification of the ordinary primes for $f(x)$. The precise classification, however, is a little bit complicated to describe.

Similarly, we can consider the Laurent polynomial

$$f(x) = a_1 x_1 + \cdots + a_n x_n + a_{n+1} x_1^{v_1} \cdots x_n^{v_n}, \quad a_j \neq 0, v_j > 0, \tag{4.4}$$

or by making the invertible change of variables $x_i \to x_i^{-1}$, we get the new equivalent Laurent polynomial

$$f(x) = a_1 \frac{1}{x_1} + \cdots + a_n \frac{1}{x_n} + a_{n+1} \frac{1}{x_1^{v_1} \cdots x_n^{v_n}}, \quad a_j \neq 0, v_j > 0. \tag{4.5}$$

The polynomial $f(x)$ in (4.4)-(4.5) is non-degenerate if and only if none of the $v_j$ is divisible by $p$. The L-function $L^*(f, T)^{(-1)^{n-1}}$ in this case is a polynomial of degree $v_1 + \cdots + v_n$. It is ordinary at $p$ if $p - 1$ is divisible by $\mathrm{lcm}(v_1, \cdots, v_n)$.

For another example with $k > 1$, we consider the sum of two polynomials from (4.4)-(4.5). Namely, let $f(x)$ be the $n$-dimensional ($n > 1$) generalized bi-Kloosterman polynomial

$$f(x) = a_1 x_1 + \cdots + a_n x_n + a_{n+1} \frac{1}{x_1^{u_1} \cdots x_n^{u_n}} + b_1 \frac{1}{x_1} + \cdots + b_n \frac{1}{x_n} + b_{n+1} x_1^{v_1} \cdots x_n^{v_n}, \tag{4.6}$$

where the coefficients are non-zero and the $u_i, v_j$ are positive integers. The exponential sum associated to this polynomial (in the special case $u_i = v_j = 1$) arises from Kim's [Ki] calculation of the Gauss sum for the unitary group. One checks that $\Delta(f)$ has exactly $(2^n - 2) + 2n = 2^n + 2n - 2$ co-dimension 1 faces $\delta_i$ not containing the origin. The largest invariant factors (in fact the determinants) of these faces are given by $u_i$, $v_j$ and 1 ($2^n - 2$ of them). Thus, the polynomial $f(x)$ in (4.6) is non-degenerate if and only if none of the $u_i, v_j$ is divisible by $p$. The L-function $L^*(f,T)^{(-1)^{n-1}}$ in this case is a polynomial of degree $2^n + 2n - 2$. Corollary 4.4 shows that $f$ is ordinary at $p$ if $p - 1$ is divisible by $\mathrm{lcm}(u_1, \cdots, u_n, v_1, \cdots, v_n)$. In particular, if all $u_i = v_j = 1$ as in [Ki], then $f$ is ordinary for every prime $p$.

We now turn to some easy applications of Theorem 4.1 to the Adolphson-Sperber conjecture. In the case $n = 1$, Corollary 4.4 with $k \leq 2$ immediately implies

**Corollary 4.5.** *Conjecture 2.9 holds for $n = 1$.*

In the case $n = 2$, if $\Delta$ has only one co-dimension 1 face not containing the origin, then $\Delta$ must be a simplex. This fact together with Theorem 4.1 gives the following weaker version of Conjecture 2.9 for $n \leq 2$.

**Corollary 4.6.** *Theorem 2.10 holds if $n \leq 2$.*

To treat the full form of Conjecture 2.9 even in the case $n = 2$ or its weaker version in higher dimension $n \geq 3$, we need to introduce further decomposition theorems. We will describe the new collapsing decomposition in next section.

### 4.2. Collapsing decomposition for generic Newton polygon.

We will assume that $n \geq 2$ as the case $n = 1$ is already handled by the facial decomposition. Let $\mathcal{V} = \{V_1, \cdots, V_J\}$ be the set of $J$ fixed lattice points in $\mathbf{R}^n$. Let $\Delta$ be the convex polyhedron in $\mathbf{R}^n$ generated by the origin and the lattice points in $\mathcal{V}$. We assume that $\Delta$ is $n$-dimensional. By the facial decomposition, we may assume that $\Delta$ has only one co-dimension 1 face $\delta$ not containing the origin and all $V_j \in \delta$. To decompose $\Delta$, we will decompose the unique face $\delta$. Actually, we will be decomposing the set $\mathcal{V}$ since we are working in a little more general setting. In the special case when $\mathcal{V}$ consists of all lattice points in $\delta$, then our decomposition can be described purely in terms of $\delta$. Clearly, the set $\mathcal{V}$ has at least $n$ elements. If the set $\mathcal{V}$ has exactly $n$ elements, then the set $\mathcal{V}$ is called **indecomposable**.

Choose an element in $\mathcal{V}$ which is a vertex of $\delta$, say $V_1$. The collapsing decomposition of $\mathcal{V}$ with respect to $V_1$ is simply the convex decomposition of $\mathcal{V}$ resulted from the collapsing after removing the vertex $V_1$. We now describe the collapsing decomposition more precisely. Let $\mathcal{V}_1 = \mathcal{V} - \{V_1\}$ be the complement of $V_1$ in $\mathcal{V}$. Let $\delta_1$ be the convex polytope of the lattice points in $\mathcal{V}_1$. This is a subset of $\delta$. Let $\delta_1'$ be the closure of $\delta - \delta_1$. Then, the intersection $\delta_1 \cap \delta_1'$ consists of finitely many different co-dimension 1 faces $\{\epsilon_2, \cdots, \epsilon_k\}$ of $\delta_1$. Let $\delta_i$ ($2 \leq i \leq k$) be the convex closure of $\epsilon_i$ and $V_1$. Then, each $\delta_i$ is $(n-1)$-dimensional. Let $\mathcal{V}_i$ ($2 \leq i \leq k$) be the intersection $\mathcal{V} \cap \delta_i$. Then, each $V_j$ lies in one (possibly more) of the subsets $\mathcal{V}_i$ of $\mathcal{V}$. The collapsing decomposition of $\mathcal{V}$ with respect to $V_1$ is defined to be

$$\mathcal{V} = \cup_{i=1}^k \mathcal{V}_i. \tag{4.7}$$

The full collapsing decomposition theorem for the degree polygon applies to the decomposition in (4.7). To describe it only in terms of Newton polygon, we must further reduce to indecomposable diagonal situation, in which case the degree polygon coincides with its upper bound if and only if the Newton polygon coincides with its lower bound. And thus, we will be able to replace the degree polygon by the simpler Newton polygon.

Applying the collapsing decomposition in (4.7) to each $\mathcal{V}_i$ by choosing a vertex of $\mathcal{V}_i$, by induction, we will eventually be able to decompose $\mathcal{V}$ as a finite union of indecomposable ones:

$$\mathcal{V} = \cup_{i=1}^{m} \mathcal{W}_i, \tag{4.8}$$

where each $\mathcal{W}_i$ has exactly $n$ elements and thus indecomposable. Furthermore, the subsets $\mathcal{W}_i$'s are different although they may have non-empty intersections. The complete decomposition in (4.8) resulting from a sequence of collapsing decompositions is called a **complete collapsing decomposition** of $\mathcal{V}$. Note that such a complete collapsing decomposition is in general not unique because of the choice of a vertex in each stage of the collapsing construction.

Let $f$ be the generic Laurent polynomial associated to the set $\mathcal{V}$:

$$f = \sum_{j=1}^{J} a_j x^{V_j}.$$

Let the $\mathcal{W}_i$ ($1 \leq i \leq m$) be a complete collapsing decomposition of $\mathcal{V}$ as in (4.8). For each integer $1 \leq i \leq m$, let $f_i$ be the restriction of $f$ to $\mathcal{W}_i$:

$$f_i = \sum_{V_j \in \mathcal{W}_i} a_j x^{V_j}.$$

This is the generic Laurent polynomial associated to the indecomposable $\mathcal{W}_i$. We have

**Theorem 4.7 (collapsing decomposition).** *Let the $\mathcal{W}_i$ ($1 \leq i \leq m$) be a complete collapsing decomposition of $\mathcal{V}$. If each $f_i$ ($1 \leq i \leq m$) is generically ordinary for some prime $p$, then $f$ is also generically ordinary for the same prime $p$.*

Since each $\mathcal{W}_i$ is indecomposable, each $f_i$ is the generic diagonal Laurent polynomial whose exponents are in $\mathcal{W}_i$. As a consequence of Theorem 4.7 and Corollary 3.5, we obtain

**Corollary 4.8.** *Let $M_i$ be the non-singular $n \times n$ matrix whose column vectors are the elements of $\mathcal{W}_i$. Let $d_{ni}$ be the largest invariant factor of the matrix $M_i$. Let*

$$D^* = \mathrm{lcm}\{d_{n1}, d_{n2}, \cdots, d_{nm}\}. \tag{4.9}$$

*If $p - 1$ is divisible by $D^*$, then $f$ is generically ordinary for the prime $p$.*

The number $D^*$ in (4.9) depends on our choice of the complete collapsing decomposition. Although any choice would yield a good upper bound for the smallest $D^*$, we do not know a general algorithm which tells us how to choose a best complete collapsing decomposition. Taking $\mathcal{V}$ to be the full set of lattice points in $\delta$, Corollary 4.8 reduces

**Corollary 4.9.** *The weaker version of the Adolphson-Sperber conjecture as stated in Theorem 2.10 is true in every dimension.*

Taking $n \leq 3$ in Corollary 4.8, we deduce from Corollary 3.7, the following full form of the Adolphson-Sperber conjecture for $n \leq 3$.

**Corollary 4.10.** *Conjecture 2.9 is true for $n \leq 3$.*

Applying this Corollary to the 3-variable polynomials of the form $x_3 f(x_1, x_2)$, we can deduce

**Corollary 4.11.** *Let $\Delta$ be a 2-dimensional integral convex polyhedron in $\mathbf{R}^2$. Then, for every prime $p$, the zeta function of the affine toric curve $f(x_1, x_2) = 0$ $(x_i \neq 0)$ is ordinary for a generic non-degenerate $f$ with $\Delta(f) = \Delta$.*

A simple boundary argument (see [W3]) then gives

**Corollary 4.12.** *Let $\Delta$ be a 2-dimensional integral convex polyhedron in the first quadrant of $\mathbf{R}^2$. Assume that $\Delta$ contains a non-zero vertex on both the $x_1$ axis and the $x_2$-axis. Then, for every prime $p$, the zeta function of the affine curve $f(x_1, x_2) = 0$ is ordinary for a generic non-degenerate $f$ with $\Delta(f) = \Delta$.*

This result generalizes the well known fact that a generic non-degenerate plane curve of degree $d$ is ordinary, see [Ko] and [Mi]. More generally, in next section, we will show that a generic non-degenerate affine hypersurface of degree $d$ is ordinary (affine version of Mazur's conjecture). But Corollary 4.11 and Corollary 4.12 do not generalize to higher dimensional $\Delta$. In fact, the 5-dimensional counter-example in Section 3.4 shows that the higher dimensional generalization of Corollary 4.11 is already false for the 3-dimensional generic toric affine hypersurface $f(x_1, \cdots, x_4)$ with $\Delta(f) = \Delta$ for some 4-dimensional $\Delta$.

### 4.3. Hyperplane decomposition for generic Newton polygon.

In this section, we describe a weaker consequence of the hyperplane decomposition theorem from [W3]. As an application, we explain how it implies Mazur's conjecture for the zeta function of a generic hypersurface. Again, with no loss of generality, we can assume that $n \geq 2$.

Let $\mathcal{V} = \{V_1, \cdots, V_J\}$ be the set of $J$ fixed lattice points in $\mathbf{R}^n$. Let $\Delta$ be the convex polyhedron in $\mathbf{R}^n$ generated by the origin and the lattice points in $\mathcal{V}$. We assume that $\Delta$ is $n$-dimensional. By the facial decomposition, we may assume that $\Delta$ has only one co-dimension 1 face $\delta$ not containing the origin and all $V_j \in \delta$. To decompose $\Delta$, we will decompose the unique face $\delta$. Actually, we will be decomposing the set $\mathcal{V}$ since we are working in a little more general setting. In the special case when $\mathcal{V}$ consists of all lattice points in $\delta$, then our decomposition can be described purely in terms of $\delta$.

Let $k$ be a positive integer. Let $H_i$ $(1 \leq i \leq k)$ be a sequence of successive parallel $(n-2)$-dimensional hyperplanes contained in the $(n-1)$-dimensional hyperplane spanned by $\delta$. We assume that the initial hyperplane $H_1$ does not intersect the interior of $\delta$. That is, the intersection of $H_1 \cap \delta$ is a (possibly empty) face of $\delta$. These hyperplanes $H_i$ cut $\delta$ into finitely many pieces. Let $\delta_i$ $(1 \leq i \leq k-1)$ be the closed piece of $\delta$ which is contained between $H_i$ and $H_{i+1}$. Let $\delta_k$ be the remaining closed piece which is the closure of $\delta - \cup_{i=1}^{k-1} \delta_i$. We assume two further conditions on the $\delta_i$. First, each $\delta_i$ ($1 \leq i \leq k$) is an $(n-1)$-dimensional **integral** polytope. Second, there are no lattice points in the interior of each $\delta_i$. Under these assumptions on $\delta_i$, we call the hyperplane decomposition

$$\delta = \cup_{i=1}^{k} \delta_i. \tag{4.10}$$

an **admissible** hyperplane decomposition of $\delta$. It should be noted that such an admissible hyperplane decomposition does not always exist.

Assume that we are given an admissible hyperplane decomposition of $\delta$ as in (4.10). Let $\mathcal{W}_i$ be the intersection $\mathcal{V} \cap \delta_i$. Assume that for each $1 \leq i \leq k$, $\delta_i$ coincides with the convex polyhedra generated by the lattice points in $\mathcal{W}_i$. Then the decomposition

$$\mathcal{V} = \cup_{i=1}^{k} \mathcal{W}_i \tag{4.11}$$

is called an **admissible** hyperplane decomposition of the set $\mathcal{V}$. Let $f$ be the generic Laurent polynomial associated to the set $\mathcal{V}$:

$$f = \sum_{j=1}^{J} a_j x^{V_j}.$$

For each integer $1 \leq i \leq k$, let $f_i$ be the restriction of $f$ to $\mathcal{W}_i$:

$$f_i = \sum_{V_j \in \mathcal{W}_i} a_j x^{V_j}.$$

This is the generic Laurent polynomial associated to the set $\mathcal{W}_i$. We have the following result which is a consequence of Theorem 7.1 in [W3].

**Theorem 4.13 (hyperplane decomposition).** *Suppose that we are given an admissible hyperplane decomposition of $\mathcal{V}$ as in (4.11). If each $f_i$ ($1 \leq i \leq k$) is generically ordinary for some prime $p$, then $f$ is also generically ordinary for the same prime $p$.*

Since each piece $f_i$ is not necessarily diagonal yet, we may need to apply the above decomposition theorems several times in order to reduce to the diagonal case.

As an example, with $n$ replaced by $n+1$, we consider the case that $\Delta$ is the $(n+1)$-dimensional simplex in $\mathbf{R}^{n+1}$ generated by the origin, $V$ and $V + de_i$ ($1 \leq i \leq n$), where $d$ is a fixed positive integer and $V$ is a fixed lattice point on the hyperplane $x_{n+1} = D$ for some positive integer $D$. The $e_i$ are the standard unit vectors in $\mathbf{R}^n$ embedded in $\mathbf{R}^{n+1}$ via the equation $x_{n+1} = 0$. Let $\mathcal{V}$ be the set of all lattice points on the face $\delta$ generated by the vertices $V$ and $V + de_i$. Under these assumptions, we have

**Theorem 4.14.** *Let $f$ be the generic Laurent polynomial whose exponents are in $\mathcal{V}$. Then $f$ is generically ordinary at prime $p$ if $p \equiv 1 \pmod{D}$.*

**Proof**. To simplify notations, we may assume that $\delta$ is the $n$-dimensional simplex in $\mathbf{R}^n$ generated by the origin and the vertices $de_i$. That is, $V = (0, \cdots, 0, D)$. Now, we use the regular subdivision of a simplex as described in [KK, p117-119]. Define the cumulative coordinates as follows:

$$y_1 = x_1,$$
$$y_2 = x_1 + x_2,$$
$$\cdots\cdots$$
$$y_n = x_1 + x_2 + \cdots + x_n.$$

One sees that $\delta$ is identified with the set of $n$-tuples of real numbers $y_1, \cdots, y_n$ such that

$$0 \leq y_1 \leq y_2 \leq \cdots \leq y_n \leq d.$$

Consider the following hyperplanes in $\mathbf{R}^n$:

$$H(i,j;k): \; y_i - y_j = k, \text{ where } 0 \leq j < i \leq n, \; 0 \leq k \leq d. \tag{4.12}$$

We claim that these hyperplanes cut $\delta$ into integral simplexes with determinant 1 (equivalently, with volume $1/n!$). To prove the claim, let $y_1, \cdots, y_n$ be a point in $\delta' = \delta - \cup_{i,j,k} H(i,j;k)$. If we put $t_i = y_i - [y_i]$, it is clear that $t_i \neq t_j$ whenever $i > j$ (if $t_i = t_j$, then $y_i - y_j$ is an integer $k$ in the interval $[0,d]$. This contradicts with the assumption that $y_1, \cdots, y_n \in \delta'$). Hence there is a unique permutation $\pi$ such that

$$0 < t_{\pi(1)} < t_{\pi(2)} < \cdots < t_{\pi(n)} < 1.$$

It is also clear that the connected component of $\delta'$ determined by the point $y_1, \cdots, y_n$ is the set of points $y'_1, \cdots, y'_n$ such that

$$0 < y'_{\pi(1)} - [y'_{\pi(1)}] < \cdots < y'_{\pi(n)} - [y'_{\pi(n)}] < 1.$$

This is clearly the interior of a simplex with integral vertices and with determinant 1. The claim is proved. Let $\delta_i$ be any integral $n$-dimensional simplex in $\delta$ with determinant 1. Let $\Delta_i$ be the simplex generated by $\delta_i$ and the origin. Then, the largest invariant factor of $\Delta_i$ is exactly $D$ since the face $\delta_i$ has determinant 1. Thus, by Corollary 3.5, the diagonal Laurent polynomial associated to $\Delta_i$ is ordinary if $p-1$ is divisible by $D$.

Now, for each fixed pair $(i,j)$ with $i > j$ as in (4.12), the hyperplanes $H(i,j;k)$ ($0 \leq k \leq d$) give arise an admissible hyperplane decomposition of $\delta$ (respectively of $\mathcal{V}$). Using the claim together with Theorem 4.13, by induction we deduce that Theorem 4.14 is true.

Theorem 4.14 gives a higher dimensional example for which the full form of the Adolphson-Sperber conjecture is true. For the application to the zeta function of a generic affine hypersurface $f(x_1, \cdots, x_n)$ of degree $d$, we need to consider the L-function of the new generic non-degenerate polynomial $x_{n+1} f(x_1, \cdots, x_n)$. Taking $D = 1$ in Theorem 4.14 and using a simple boundary argument ([W3]), we conclude that the the zeta function of a generic non-degenerate hypersurface of degree $d$ is ordinary for every prime $p$. This is the affine version of Mazur's conjecture.

For the application to the zeta function of a generic projective hypersurface of degree $d$, we need to consider the case that $\Delta$ is the $n$-dimensional simplex in $\mathbf{R}^{n+1}$ generated by the origin and $V + de_i$. This $\Delta$ is a co-dimension 1 face of the polyhedron treated in Theorem 4.14. The essentially trivial boundary decomposition theorem (Corollary 5.2 in [W3]) shows that on the chain level, we still have the generic ordinary property if $p-1$ is divisible by $D$. Although $\Delta$ does not have the maximal dimension $n+1$, Adolphson-Sperber's results (Theorem 2.1 and Theorem 2.8) carry over to this case. We restricted our description to the maximal dimensional case so that the notations are simpler. Thus, the chain level result passes to the L-function. We conclude that the L-function is generically ordinary if $p-1$ is divisible by $D$. Taking $D=1$, one gets Mazur's original conjecture that the zeta function of a generic projective hypersurface of degree $d$ is ordinary for every prime $p$, see [W3] for more details. Actually, if one wishes, one can also use induction to derive the projective version directly from the above affine version.

One can construct many more higher dimensional examples for which the full Adolphson-Sperber conjecture is true. For instance, if $d_1, \cdots, d_n$ are $n$-positive integers, let $f(x_1, \cdots, x_n)$ be the generic polynomial which has degree $d_i$ with respect to the $i$-th variable $x_i$ for each $1 \leq i \leq d$. Then the zeta function of the generic non-degenerate affine hypersurface defined by $f$ is ordinary for every prime $p$. In this case, the face $\delta$ is an integral box and $D=1$. One can first use the hyperplane decomposition to reduce to the case of a unit box. One can then apply the star decomposition as did in [W3] or the new collapsing decomposition as described in the previous section.